\documentclass[a4paper]{article}

\usepackage{amsfonts}
\usepackage{amsmath}
\usepackage{ntheorem}
\usepackage{hyperref}

\newtheorem{conj}{Conjecture}[section]
\newtheorem{theo}[conj]{Theorem}
\newtheorem{prop}[conj]{Proposition}

\begin{document}

\title{Note on a Conjecture of Wegner}
\author{Dominik Kenn}
\date{February 2, 2008} %Datum der erstmaligen Veröffentlichung bei ArXiv
\maketitle

\begin{abstract}
The optimal packings of $n$ unit discs in the plane are known for those $n \in \mathbb N$, which satisfy certain number theoretic conditions. Their geometric realizations are the \textit{extremal Groemer packings} (or \textit{Wegner packings}). But an extremal Groemer packing of $n$ unit discs does not exist for all $n\in\mathbb N$ and in this case, the number $n$ is called \textit{exceptional}. \\
We are interested in number theoretic characterizations of the exceptional numbers. \\
A counterexample is given to a conjecture of Wegner concerning such a characterization. We further give a characterization of the exceptional numbers, whose shape is closely related to that of Wegner's conjecture.
\end{abstract}
\section{Introduction}
It has been conjectured by L. F. T\'oth \cite{Tot}, that for any given $n\in\mathbb N$, the optimal packing with respect to the convex hull of $n$ discs has to be a subset with $n$ elements of the hexagonal lattice packing whose convex hull is "as similar to a regular hexagon as possible". In this context it is clear, that a packing is not optimal, if a further disc fits into the convex hull of the previous discs without intersecting the interior of at least one of the previous discs. Therefore we will restrict ourselves to \textit{dense} packings. Beyond this, it suffices to consider only unit discs.\\
Wegner \cite{Weg1} solved the above problem almost completely. But the packing of $n$ unit discs minimizing the area of the convex hull of the discs is not known for all $n\in\mathbb N$. Wegner proved that, in case of existence, the so-called \textit{extremal Groemer packings} are optimal and they are the only optimal packings. For $n$ with $n=1+6\binom a2$, the extremal Groemer packing of $n$ unit discs exists and indeed has the shape of a regular hexagon and there is no other packing with the same density. (This has been part of T\'oth's conjecture.) If there exists an extremal Groemer packing for some $n$ with $n=1+6\binom a2+ab+c$, then it has the shape of a degenerated hexagon. In this case, there can be more than one extremal Groemer packing of $n$ unit discs. \\
Now a number $n$ is called \textit{exceptional}, if there is no extremal Groemer packing of $n$ unit discs. Since the existence of extremal Groemer packings can be formulated as a number theoretic problem, it is natural to ask for a number theoretic characterization of the exceptional numbers. Moreover, the knowledge of such a characterization may be useful to describe the geometric shape of the optimal packings of $n$ unit discs for exceptional numbers $n$. \\
Wegner \cite{Weg2} gave a conjecture concerning a number theoretic characterization of the exceptional numbers. B\"or\"oczky and Ruzsa \cite{BoeR} proved another characterization, but surprisingly, they did not compare their results to the results of Wegner. \\
We prove here that Wegner's conjecture is wrong and by means of investigating the connection between the results of Wegner and B\"or\"oczky and Ruzsa, we will "correct" Wegner's conjecture.
\section{Notations and Preliminaries}
We recall some necessary facts and definitions. \\
Each $n \in \mathbb N$ can be written as
\[
  n=1+6 \binom a2 +ab+c
\]
where the parameters $a,b,c \in \mathbb Z$ have to be chosen maximal in this order. It then follows that
\[
  1\le a, \qquad 0\le b\le 5 \qquad \text{and} \qquad 0\le c<a \,.
\]
Set
\[
  p_0(n):=6(a-1)+b+1-\delta_{0,b+c}
\]
where $\delta$ denotes the Kronecker symbol.
\\[3ex]
A packing of $n\in\mathbb N$ unit discs $B^2$ in Euclidian space $\mathbb R^2$ is denoted by $C_n+B^2$, where $C_n$ is the set of the centers of the discs. A packing $C_n+B^2$ is called a \textit{Groemer packing}, if the elements of $C_n$ form a, not necessarily regular, hexagon in $\mathbb R^2$. Let $conv(M)$ denote the convex hull of a set $M$ and $\partial M$ its boundary. For a Groemer packing $C_n+B^2$, a disc $B^2 \subset C_n+B^2$ is called a \textit{boundary disc}, if $\partial B^2\ \cap\ \partial (conv(C_n+B^2)) \ne \emptyset$. Now let $p(C_n+B^2)$ denote the number of boundary discs of $C_n+B^2$. Wegner proved $p(C_n+B^2) \ge p_0(n)$ for each Groemer packing. If equality holds, the Groemer packing is called \textit{extremal} (or a \textit{Wegner packing}). Combining the results of \cite{Weg3} and \cite{Har}, one gets
\begin{equation} \label{eq_p0n}
  p_0(n)=\Bigl\lceil \sqrt{12n-3} \Bigr\rceil -3 \,.
\end{equation}
For $n$ not too small, each set $conv(C_n)$ related to a Groemer packing $C_n+B^2$ has six sides, that means there are six line segments containing the centers of the boundary discs. For each of those line segments, let $p_i,\ 1\le i\le 6,$ denote the number of centers lying on that line segment. Then one gets a sequence $p_1,\dots,p_6$, called the \textit{boundary sequence} of the packing. The boundary sequence is uniquely determined by the packing up to cyclic permutations. \\
Now $p_i+p_{i+1}=p_{i+3}+p_{i+4}$ and a simple consideration shows
\begin{equation} \label{eq_n}
  n=(p_1+p_2-1)(p_3+p_4-1)-\binom{p_1}2-\binom{p_4}2
\end{equation}
and for extremal Groemer packings, one also has
\begin{equation} \label{eq_p0}
  p_0(n)=p_1+2p_2+2p_3+p_4-6 \,.
\end{equation}
For each Groemer packing $C_n+B^2$ with $n=1+6\binom a2+ab+c$ and $b$ and $c$ not both equal $0$, the numbers $p_1,\dots,p_6$ have to satisfy the following conditions (see \cite{Weg1}, page 6):
\begin{equation} \label{eq_pis}
  \frac{a-1}2 \le p_i \le 2a-c \,.
\end{equation}
So an extremal Groemer packing $C_n+B^2$ exists if and only if there are natural numbers $p_1,\dots,p_4$ with (\ref{eq_n}), (\ref{eq_p0}) and (\ref{eq_pis}). Those $n\in \mathbb N$, for which no extremal Groemer packing exists, are called \textit{exceptional numbers}. In \cite{Weg2}, Wegner gives the following
\begin{conj} \label{conj_wegner}
  \textbf{\textup{(Wegner)}} A number $n\in\mathbb N$ is exceptional if and only if the parameters $a,b,c$ satisfy one of the following conditions:
  \begin{align*}
    i)\quad & b=2 \quad\text{and}\quad a-c\equiv-6m \mod 9^{m+1} \text{ for } m \in \mathbb N_0 \\
   ii)\quad & b=5 \quad\text{and}\quad a-c\equiv6\cdot 9^m \mod 9^{m+1} \text{ for } m \in \mathbb N_0 \,. 
  \end{align*} 
\end{conj}
In \cite{BoeR}, B\"or\"oczky and Ruzsa prove the following
\begin{theo} \label{theo_boer}
  \textbf{\textup{(B\"or\"oczky and Ruzsa)}} A number $n\in\mathbb N$ is exceptional if and only if
  \[
    \Bigl\lceil \sqrt{12n-3} \Bigr\rceil ^2 +3-12n = (3k-1)9^\ell
  \]
  for some $k,\ell \in \mathbb N$.
\end{theo}
\section{A counterexample to the Wegner conjecture}
\begin{prop} \label{prop_gegnbsp}
  The Wegner conjecture is wrong in general.
\end{prop}
\textbf{Proof} \\
We refer to part i) of conjecture \ref{conj_wegner}. Let $b=2$ and choose $m=2$. Then we have
\[
   a-c \equiv -12 \equiv 717 \mod 9^3 \,.
\]
Now for $a=717$ and $c=0$ we obtain
\[
  n=1+6\binom a2+ab+c=1541551 
\]
and condition i) claims that this is an axceptional number. \\
By computer-aided calculation, we find (beside other solutions)
\[
  p_1=702, \qquad p_2=717, \qquad p_3=714 \qquad \text{and} \qquad p_4=741 \,.
\]
The $p_i$'s satisfy (\ref{eq_n}), (\ref{eq_p0}) and (\ref{eq_pis}), so $n=1541551$ is no exceptional number.
\hfill$\square$ \\
[3ex]
It can be checked by an easy calculation, that there are no $k,\ell\in\mathbb N$ so that the equation from Theorem \ref{theo_boer} holds for $n=1541551$, so this number is not exceptional in the sense of B\"or\"oczky and Ruzsa. In particular, conjecture \ref{conj_wegner} and theorem \ref{theo_boer} are not equivalent.\\
Now we give a characterization of the exceptional numbers similar to the Wegner conjecture by using the result of B\"or\"oczky and Ruzsa. \\
[3ex]
For the sake of notation, let $\sum_{i=j}^k a_i =0$ for $k<j\in\mathbb Z$.
\begin{prop} \label{prop_ausn}
  A number $n\in\mathbb N$ is exceptional if and only if the parameters $a,b,c$ satisfy one of the following conditions:
  \begin{align*}
    i) \quad & b=2 \text{ and } a-c\equiv -6m \mod 9^\ell \,, \\
             & \text{where } m=\sum_{i=0}^{\ell-2}9^i \\
             & \text{for } \ell \text{ with } 9^\ell(3k-1)=12(a-c)-9 \,, \\
   ii) \quad & b=5 \text{ and } a-c \equiv 6\cdot 9^{\ell-1} \mod 9^\ell \,, \\
             & \text{for } \ell \text{ with } 9^\ell(3k-1)=12(a-c) \,.
  \end{align*}
\end{prop}
\textbf{Proof} \\
The numbers $n=1+6\binom a2$ are not exceptional, so let $b$ and $c$ be not both equal to $0$. \\
From Theorem \ref{theo_boer} it follows with (\ref{eq_p0n}) that $n\in\mathbb N$ is exceptional iff
\begin{equation} \label{eq_ausn}
  9^\ell(3k-1) = 12(a-c)+(b-2)^2-9 
\end{equation}
for some $k,\ell\in\mathbb N$. \\
[3ex]
\textbf{I)} Let $n$ be an exceptional number. Equation (\ref{eq_ausn}) modulo 3 yields $b=2$ or $b=5$, since $0\le b\le5$. \\
[3ex]
\textbf{i)} Let $b=2$. Then
\[
  4(a-c)=9^\ell k-3(9^{\ell-1}-1) \,,
\]
therefore $k=4z$ for some $z\in\mathbb N$. Hence a direct computation shows
\[
  a-c=9^\ell z-6\sum_{i=0}^{\ell-2}9^i \,.
\]
Now let $m:=\sum_{i=0}^{\ell-2}9^i$.\\
[3ex]
\textbf{ii)} Let $b=5$. Then
\[
  12(a-c)=9^\ell(3k-1)\,,
\]
hence $3k-1\equiv0 \mod 4$ and so $k\equiv3 \mod 4$. We write $k=4z+3$ for some $z\in\mathbb N_0$. That yields
\[
  a-c=9^\ell z+6\cdot9^{\ell-1} \,.
\]
\textbf{II)} Let $a,b,c$ satisfy condition i) or ii). We show that they yield exceptional numbers.  \\
[3ex]
\textbf{i)} Let $b=2$ and $a-c \equiv -6m \mod 9^{\ell}$. \\
For $\ell,z\in\mathbb N$, there is a $k\in\mathbb N$ with
\[
  -8\sum_{i=0}^{\ell-2}9^i+12\cdot9^{\ell-1}z-1=9^{\ell-1}(3k-1)\,.
\]
This is easily seen by induction. \\
We write $a-c=-6m+9^\ell z$ for some $z\in\mathbb N$. Then
\begin{align*}
  (p_0(n)+3)^2+3-12n 
    & = 9(-8\sum_{i=0}^{\ell-2}9^i+12\cdot9^{\ell-1}z-1) \\
    & = 9^\ell(3k-1) \,.
\end{align*}
This is equation (\ref{eq_ausn}) with $b=2$. \\
[3ex]
\textbf{ii)} Let $b=5$ and $a-c \equiv 6\cdot9^{\ell-1} \mod 9^\ell$. We write $a-c=6\cdot9^{\ell-1}+9^\ell z$ for some $z\in\mathbb N_0$ and get
\[
  (p_0(n)+3)^2+3-12n = 9^\ell(3(4z+3)-1) \,.
\]
This is equation (\ref{eq_ausn}) with $b=5$ and $k:=4z+3$.
\hfill$\square$
\section*{Acknowledgement}
I am grateful to Prof. A. Sarti for many helpfull discussions.
\\[3ex]
The final publication is available at \href{http://springerlink.com/}{springerlink.com}, \\
DOI 10.1007/s13366-011-0004-3.
\bibliography{biblio}

\begin{thebibliography}{1}

\bibitem{BoeR}
K.~J. B{\"o}r{\"o}czky and I.~Z. Ruzsa.
\newblock Note on an {I}nequality of {W}egner.
\newblock {\em Discrete {C}omput. {G}eom.}, 37:245--249, 2007.

\bibitem{Har}
H.~Harborth.
\newblock L{\"o}sung zum {P}roblem 664 {A}.
\newblock {\em Elem. {M}ath.}, 29:14--15, 1974.

\bibitem{Tot}
L.~F. T\'oth.
\newblock Research {P}roblem 13.
\newblock {\em Period. {M}ath. {H}ungar.}, 6:197--199, 1975.

\bibitem{Weg3}
G.~Wegner.
\newblock Zur {K}ombinatorik von {K}reispackungen.
\newblock {\em Koll. {D}iskr. {G}eom.}, 2:225--230, 1980.

\bibitem{Weg2}
G.~Wegner.
\newblock Extremale {G}roemerpackungen.
\newblock {\em Studia Sci. Math. Hungar.}, 19:299--302, 1984.

\bibitem{Weg1}
G.~Wegner.
\newblock {\"U}ber endliche {K}reispackungen in der {E}bene.
\newblock {\em Studia Sci. Math. Hungar.}, 21:1--28, 1986.

\end{thebibliography}
\bibliographystyle{plain}
Dominik Kenn, Institut f\"ur Mathematik, Johannes-Gutenberg-Universit\"at Mainz, Staudingerweg 9, 55099 Mainz \\
dominik.kenn@googlemail.com
\end{document}